\documentclass[leqno,dvips,hdvips,11pt]{amsart}

\usepackage{booktabs}
\usepackage{amsmath, amssymb, url,a4wide}
\usepackage{xcolor,graphicx}
\usepackage[]{xr-hyper,hyperref}

\newtheorem{theorem}                 {Theorem}      [section]
\newtheorem{remark}[theorem]{Remark}

\title[Minimal translation surfaces in $\mathrm{Nil}_3$]{Minimal translation surfaces in \\the Heisenberg group $\mathrm{Nil}_3$}

\author{Jun-ichi Inoguchi}
\author{Rafael L\'opez}
\author{Marian Ioan Munteanu}

\address[J.-I.~Inoguchi]
{Department of Mathematical Sciences\\
Faculty of Science\\
Kojirakawa-machi 1-4-12\\
Yamagata 990-8560, Japan
}
\email{inoguchi (at) sci.kj.yamagata-u.ac.jp}

\address[R. L\'opez]{Departamento de Geometr\'{\i}a y Topolog\'{\i}a\\
Universidad de Granada, 18071 Granada, Spain}
\email{rcamino (at) ugr.es}

\address[M.~I.~Munteanu]
{Michigan State University\\
Department of Mathematics\\
Wells Hall\\
48824-1029 East Lansing\\ USA}

\address[]
{Al.I. Cuza University of Iasi\\
Faculty of Mathematics\\
Bd. Carol I, n. 11\\
700506 Iasi\\
Romania}
\email{marian.ioan.munteanu (at) gmail.com}

\date{\today}



\begin{document}

\begin{abstract}
A translation surface in the Heisenberg group $\mathrm{Nil}_3$ is a surface constructed by multiplying
(using the group operation) two curves.
We completely classify minimal translation surfaces in the Heisenberg group $\mathrm{Nil}_3$.

\end{abstract}
\maketitle

\section{Introduction}

A surface $M$ in the Euclidean space is called a \emph{translation surface} if it is given by the graph $z(x,y)=f(x)+g(y)$,
where  $f$ and $g$ are smooth functions on some interval of ${\mathbb{R}}$. In \cite{sc},
Scherk proved that, besides the planes, the only minimal translation surfaces are  given by
$$
    z(x,y)=\frac{1}{a}\log\Big|\frac{\cos(ax)}{\cos(ay)}\Big|,
$$
where $a$ is a non-zero constant. These surfaces are now referred as
\textit{Scherk's minimal surfaces}.

The study of translation surfaces in the Euclidean space was
extended when the second fundamental form was considered as a metric on a non-developable surface.
A classification is given for translation surfaces for which the second Gaussian curvature and the mean curvature are proportional \cite{mn}.
When the ambient is the Minkowski 3-space, translation surfaces of Weingarten type are classified \cite{dgv}.
In \cite{gv}, translation surfaces with vanishing second Gaussian curvature in Euclidean and Minkowski 3-space are studied.

In the last decade,
there has been an intensive effort to
develop the theory of surfaces in homogeneous Riemannian 3-spaces
of non-constant curvature. Since the discovery of
holomorphic quadratic differential
(called generalized Hopf differential or
Abresch-Rosenberg differential) for constant mean
curvature surfaces in $3$-dimensional
homogeneous Riemannian spaces
with $4$-dimensional isometry group,
global geometry of constant mean curvature
surfaces in such spaces has been extensively studied.
We refer the survey \cite{ar},\cite{fm} or
lecture notes \cite{dhm} and references therein.
In particular, integral representation formulae for minimal
surfaces in the Heisenberg group $\mathrm{Nil}_3$ were
obtained independently in \cite{d,f,i3}.
Some fundamental examples of minimal surfaces are constructed in \cite{ikos}.
Berdinski{\u\i} and Ta{\u\i}manov \cite{bt}
gave a representation formula for
minimal surfaces in $3$-dimensional
Lie groups in terms of spinors and Dirac operators.

Translation surfaces can be defined in any
3-dimensional Lie groups equipped with left
invariant Riemannian metric. Some examples of minimal
translation surfaces in the Heisenberg group
$\mathrm{Nil}_3$ are obtained in \cite{ikos}.
In our previous paper \cite{lm}, we have classified
minimal translation surfaces in the model space $\mathrm{Sol}_3$
of solvgeometry in the sense of Thurston \cite{th}.

The purpose of this article
is to study and classify
minimal translation surfaces of $\mathrm{Nil}_3$.

\section{Preliminaries}

{\bf The Heisenberg group $\mathrm{Nil}_3$} is defined as ${\mathbb{R}}^3$ with the group operation
$$
 (x,y,z)\ast(\overline{x},\overline{y},\overline{z})
   =
    \left(x+\overline{x},\ y+\overline{y},\ z+\overline{z}+\frac{x\overline{y}}{2}-\frac{\overline{x}y}{2}\right).
$$
The identity of the group is $(0,0,0)$ and the inverse of $(x,y,z)$ is given by $(-x,-y,-z)$.
It is simply connected and connected 2-step nilpotent Lie group (\cite{th}).
The following metric is left invariant
\begin{equation}
\label{metricNil}
\tilde g=dx^2+dy^2+\Big(dz+\frac12(y\,dx-x\,dy)\Big)^2.
\end{equation}
The resulting Riemannian manifold $(\mathrm{Nil}_3,\tilde{g})$
is the model space of nilgeometry in the sense of Thurston \cite{th}.

The following vector fields form a left invariant orthonormal frame
on $\mathrm{Nil}_3$:
$$
e_1 = \partial_x-\frac y2~\partial_z,\quad
e_2 = \partial_y+\frac x2~\partial_z,\quad
e_3 = \partial_z.
$$
The geometry of $\mathrm{Nil}_3$ can be described in terms of this frame as
follows.
\begin{itemize}
\item[(i)] These vector fields satisfy the commutation relations
$$
    [e_1,e_2]= e_3, \qquad [e_2,e_3]=0, \qquad [e_3,e_1]=0.
$$
\item[(ii)] The Levi-Civita connection $\widetilde\nabla$ of $\mathrm{Nil}_3$ is given by
$$\begin{array}{lll} \widetilde{\nabla}_{e_1}e_1=0, & \widetilde{\nabla}_{e_1}e_2=\frac 12 e_3, & \widetilde{\nabla}_{e_1}e_3=-\frac 12 e_2, \\
                     \widetilde{\nabla}_{e_2}e_1=-\frac 12 e_3, & \widetilde{\nabla}_{e_2}e_2=0, & \widetilde{\nabla}_{e_2}e_3=\frac 12 e_1, \\
                     \widetilde{\nabla}_{e_3}e_1=-\frac 12 e_2, & \widetilde{\nabla}_{e_3}e_2=\frac 12 e_1, & \widetilde{\nabla}_{e_3}e_3=0.
\end{array}
$$
\item[(iii)] The Riemann-Christoffel curvature tensor $\widetilde R$ of
$\mathrm{Nil}_3$ is determined by
\begin{eqnarray*}
\widetilde{R}(X,Y)Z &=& -\frac34~(\tilde g( Y,Z) X-\tilde g( X,Z) Y)\\
&& +\Big(\tilde g( Y,e_3) \tilde g( Z,e_3) X - \tilde g( X,e_3) \tilde g(Z,e_3) Y\\
&& +\tilde g( X,e_3) \tilde g(Y,Z) e_3 - \tilde g(Y,e_3)\tilde g(X,Z) e_3\Big),
\end{eqnarray*}
for $p\in \mathrm{Nil}_3$ and $X,Y,Z\in T_p \mathrm{Nil}_3$.
\end{itemize}

The Heisenberg group $\mathrm{Nil}_3$ has a four-dimensional isometry group generated by left translations $g\mapsto hg$, $h\in \mathrm{Nil}_3$,
and by rotations around $z$-axis.
Moreover the action of the isometry group
is transitive. The Heisenberg group is represented by
$\mathrm{Nil}_3=\mathrm{Nil}_3\times \mathrm{SO}(2)/\mathrm{SO}(2)$ as
a homogeneous Riemannian $3$-space.

\medskip

{\bf Minimal surfaces in $\mathrm{Nil}_3$.}
Let $r:M\longrightarrow \mathrm{Nil}_3$ be an orientable surface, isometrically immersed in $\mathrm{Nil}_3$.
Denote by $g:=r^*\tilde g$ (resp. $\nabla$)
the induced metric (resp. Levi-Civita connection) on $M$. For later use we write down the Gauss and the Weingarten formulae
$$
\begin{array}{l}
{\mathbf{(G)}}\qquad \widetilde \nabla_XY=\nabla_XY+\sigma(X,Y)N,\quad \sigma(X,Y)=\widetilde g(\widetilde\nabla_XY,N),\\[2mm]
{\mathbf{(W)}}\qquad \widetilde \nabla_XN=-AX,
\end{array}
$$
where $X,Y$ are tangent to $M$ and $N$ is a unit normal to $M$. Here $\sigma$ denotes the scalar valued
\emph{second fundamental form} of the immersion, $A$ is known as the \emph{shape operator} associated to $N$.
The shape operator $A$ is self-adjoint with respect to the metric $g$ on $M$ and it is related to $\sigma$ by
$\sigma(X,Y)=g(AX,Y)$. The \emph{mean curvature} of the immersion is defined as $H=\frac 12~{\rm{trace}}(A)$, in any point of $M$.
At each tangent plane $T_pM$ we take a basis $\{r_u,r_v\}$, where $u,v$ are local coordinates on $M$.
Denote by $E$, $F$, $G$ the coefficients of the first fundamental form on $M$:
$$
E=\widetilde g(r_u,r_u),\quad
F=\widetilde g(r_u,r_v),\quad
G=\widetilde g(r_v,r_v).
$$
We find
$$H=\frac{G\langle N,\widetilde{\nabla}_{E_1}E_1\rangle-2F\langle N,
   \widetilde{\nabla}_{E_1}E_2\rangle+E\langle N,\widetilde{\nabla}_{E_2}E_2\rangle}{2(EG-F^2)},
$$
where $\{E_1,E_2\}$ form an arbitrary basis on the surface $M$.
As we are interested in \textit{minimal surfaces}, i.e. surfaces with $H=0$,
in the above expression of $H$, we may change $N$ by other proportional
vector $\overline{N}$. Then $M$ is a minimal surface if and only if

\begin{equation}
\label{min_cond}
G\langle \overline{N},\widetilde{\nabla}_{E_1}E_1\rangle-
2F\langle \overline{N},\widetilde{\nabla}_{E_1}E_2\rangle+
E\langle \overline{N},\widetilde{\nabla}_{E_2}E_2\rangle=0.
\end{equation}

\section{Minimal translation surfaces}

In this section we define and study translation surfaces in $\mathrm{Nil}_3$.
A \emph{translation surface} $M(\gamma_1,\gamma_2)$ in $\mathrm{Nil}_3$ is a surface parametrized by
$$r:M\longrightarrow \mathrm{Nil}_3\ ,\ r(x,y)=\gamma_1(x)*\gamma_2(y),
$$
where $\gamma_1$ and $\gamma_2$ are curves situated in the planes of coordinates of ${\mathbb{R}}^3$.
Since the multiplication $*$ is not commutative, for each choice of curves $\gamma_1$ and $\gamma_2$
we may construct two translation surfaces, namely $M(\gamma_1, \gamma_2)$ and $M(\gamma_2, \gamma_1)$,
which are different. Consequently, we distinguish six types of translation surfaces in $\mathrm{Nil}_3$.

\subsection{Surfaces of type 1}

Let the curves $\gamma_1$ and $\gamma_2$ be given by $\gamma_1(x)=(x,0,u(x))$ and $\gamma_2(y)=(0,y,v(y))$.
The translation surface $M(\gamma_1,\gamma_2)=\gamma_1(x)*\gamma_2(y)$ of type 1 is thus parameterized as
$$
r(x,y)=(x,0,u(x))*(0,y,v(y))=\Big(x,y,u(x)+v(y)+\frac{xy}2\Big).
$$

The minimality condition \eqref{min_cond} yields the following ODE
\begin{equation}
\label{min_cas1:eq}
u''(x)\big(1+v'(y)^2\big)-\big(u'(x)+y\big)v'(y)+v''(y)\big[1+\big(u'(x)+y\big)^2\big]=0.
\end{equation}
In order to solve it, divide first by $1+v'(y)^2\neq0$. We get
$$
u''(x)-\big(u'(x)+y\big)~\frac{v'(y)}{1+v'(y)^2}+\frac{v''(y)}{1+v'(y)^2}~\big[1+\big(u'(x)+y\big)^2\big]=0,
$$
for all $x,y$ in the domain or $r$. Taking the derivative with respect to $x$, we obtain
$$
u'''(x)-u''(x)~\frac{v'(y)}{1+v'(y)^2}+2~\frac{v''(y)}{1+v'(y)^2}~\big(u'(x)+y\big)u''(x)=0.
$$
The case $u''(x)=0$ will be treated separately. Let us suppose now that $u''(x)\neq0$ on an open interval, and divide by $u''(x)$.
It follows
$$
\frac{u'''(x)}{u''(x)}-\frac{v'(y)}{1+v'(y)^2}+2~\frac{v''(y)}{1+v'(y)^2}~\big(u'(x)+y\big)=0.
$$

Taking the derivative with respect to $y$, we get

\begin{equation}
\label{2:eq}
-\frac {{\rm d}}{{\rm d} y}\left(\frac{v'(y)}{1+v'(y)^2}\right)+2\big(u'(x)+y\big)~\frac {{\rm d}}{{\rm d} y}\left(\frac{v''(y)}{1+v'(y)^2}\right)
    +2~\frac{v''(y)}{1+v'(y)^2}=0.
\end{equation}
If $\frac {{\rm d}}{{\rm d} y}\left(\frac{v''(y)}{1+v'(y)^2}\right)\neq0$ it follows that $u'(x)+y$ depends only on
$y$ yielding $u'(x)=$constant, and hence
$u''(x)=0$, which is a contradiction.

So $\frac{v''(y)}{1+v'(y)^2}=A$, $A\in{\mathbb{R}}$.
Replacing it in \eqref{2:eq}, one obtains
$$
-\frac{v''(y)}{1+v'(y)^2}+\frac{2v'(y)^2v''(y)}{(1+v'(y)^2)^2}+2A=0,
$$
equivalently to
$$
A\left(1+\frac{2v'(y)^2}{1+v'(y)^2}\right)=0.
$$
It follows that $A=0$, and hence $v''(y)=0$.

Summarizing, we proved that for a minimal translation surface of type 1, we should have either $u''=0$ or $v''=0$.

Moreover, if $v''=0$, then $v'(y)=c$, $c\in{\mathbb{R}}$. Replacing in \eqref{min_cas1:eq}, one obtains
$$
u''(x)(1+c^2)-c\big(u(x)+y\big)=0\ ,\ \forall\ x,y.
$$
Hence, $c=0$ and $u''(x)=0$.
Consequently, for any minimal translation surface of type 1, we have $u''(x)=0$.

Take $u(x)=ax+u_0$, $a,u_0\in{\mathbb{R}}$. Replacing in \eqref{min_cas1:eq} we get
$$
-(a+y)v'(y)+v''(y)\left(1+(a+y)^2\right)=0,
$$
which has the solution
\begin{equation}
\label{sol_v:eq}
v(y)=c\left[(a+y)\sqrt{1+(a+y)^2}+\ln(a+y+\sqrt{1+(a+y)^2})\right]+v_0\ , \ c,v_0\in{\mathbb{R}}.
\end{equation}
We conclude with the following
\begin{theorem}
\label{th1}
Minimal translation surfaces of type $1$ in the Heisenberg group $\mathrm{Nil}_3$ are parameterized by
\begin{equation}
\label{trans1:param}
r(x,y)=(x,0,u(x))*(0,y,v(y))=\left(x,y,u(x)+v(y)+\frac{xy}2\right),
\end{equation}
where $u(x)=ax+u_0$, $a,u_0\in{\mathbb{R}}$ and $v(y)$ is given by \eqref{sol_v:eq}.
\end{theorem}

\subsection{Surfaces of type 4}
Considering the same curves as in previous case, define the translation surface
of type 4 by $M(\gamma_2,\gamma_1)=\gamma_2(y)*\gamma_1(x)$. Similar computations
yield the following equation
$$
\big(1+u'(x)^2\big)v''(y)+\big(v'(y)-x\big)u'(x)+\big(1+(v'(y)-x)^2\big)u''(x)=0.
$$
We may state the next result.
\begin{theorem}
\label{th4}
Minimal translation surfaces of type $4$ in the Heisenberg group $\mathrm{Nil}_3$ are parameterized by
\begin{equation}
\label{trans4:param}
r(x,y)=(0,y,v(y))*(x,0,u(x))=\Big(x,y,u(x)+v(y)-\frac{xy}2\Big).
\end{equation}
where
$$u(x)=c\left[(a+x)\sqrt{1+(a+x)^2}+\ln(a+x+\sqrt{1+(a+x)^2})\right]+u_0,$$
and $v(y)=-ay+v_0$, with $a,c,u_0, v_0\in{\mathbb{R}}$.
\end{theorem}

\begin{remark}\rm
If $a$ and $u_0$ (resp. $a$ and $v_0$) vanish in Theorem~\ref{th1} (resp. Theorem~\ref{th4}),
then the surface $M$ is a \emph{left cylinder} over the curve $\gamma_2(y)=(0,y,v(y))$
(resp. over the curve $\gamma_1(x)=(x,0,u(x))$). See \cite{ikos}.
\end{remark}

\subsection{Surfaces of type 2}
Let us take the two curves as follows:
$\gamma_1(x)=(x,0,u(x))$ and $\gamma_2(y)=(v(y),y,0)$.
Then, the translation surface $M(\gamma_1,\gamma_2)$ of type 2 is parametrized by
$$
r(x,y)=\gamma_1(x)*\gamma_2(y)=\Big(x+v(y),y,u(x)+\frac{x y}2\Big).
$$

\smallskip

The minimality condition \eqref{min_cond} becomes:
$$
\begin{array}{rcl}
2v(y)\Big(y+u'(x)-y v'(y)u''(x)\Big)-2v'(y)\Big(2+y^2+y u'(x)\Big)&&\\[2mm]
-2v''(y)\Big(y+2u'(x)\Big)\Big(1+y^2+2yu'(x)+u'(x)^2\Big)&&\\[2mm]
+u''(x)\Big(4+v(y)^2+(4+y^2)v'(y)^2\Big)&=&0.
\end{array}
$$
Let us make some notations:
\begin{equation}
\label{eq:not_T}
\begin{array}{l}
T_0(y)=4+4v'(y)^2+\big(v(y)-y v'(y))^2,\\[2mm]
T_1(y)=2\big[v(y)-y v'(y)-2(1+2y^2)v''(y)\big],\\[2mm]
T_2(y)=-10y v''(y),\\[2mm]
T_3(y)=-4v''(y),\\[2mm]
T_4(y)=2y v(y)-2(2+y^2)v'(y)-2y(1+y^2)v''(y).
\end{array}
\end{equation}
The previous equation may be rewritten as
\begin{equation}
\label{min_cas2:eq}
T_0(y)u''(x)+T_1(y)u'(x)+T_2(y)u'(x)^2+T_3(y)u'(x)^3+T_4(y)=0.
\end{equation}
Since $T_0$ cannot vanish, divide \eqref{min_cas2:eq} by $T_0$
and then take the derivatives with respect to $x$ and $y$ respectively.
We get
$$
\frac{{\rm d}}{{\rm d} y}\left(\frac{T_1}{T_0}\right) u''(x)+
2\frac{{\rm d}}{{\rm d} y}\left(\frac{T_2}{T_0}\right) u'(x)u''(x)+
3\frac{{\rm d}}{{\rm d} y}\left(\frac{T_3}{T_0}\right) u'(x)^2u''(x)=0.
$$
The case $u''(x)=0$ will be discussed separately.
Suppose now that $u''(x)\neq0$. Dividing by $u''(x)$,
and continuing the procedure we obtain
$$
\frac{{\rm d}}{{\rm d} y}\left(\frac{T_1}{T_0}\right)=0,\quad
\frac{{\rm d}}{{\rm d} y}\left(\frac{T_2}{T_0}\right)=0, \quad
\frac{{\rm d}}{{\rm d} y}\left(\frac{T_3}{T_0}\right)=0.
$$
Looking at the form of $T_2$ and $T_3$ we conclude that $v''(y)=0$,
hence $v(y)=a y+d$, with $a,d\in{\mathbb{R}}$. The minimality condition
becomes
$$
-4a+2dy+2du'(x)+(4+4a^2+d^2)u''(x)=0,
$$
for all $x$ and $y$. It follows that $d=0$ and $u(x)=\frac{a}{2(1+a^2)}~x^2+bx+c$,
with $b,c\in{\mathbb{R}}$.

Return to the remained case $u''(x)=0$. Then $u(x)=ax+u_0$ and the minimality condition becomes
$$
(a+y)v(y)-[2+y(a+y)]v'(y)-(2a+y)[(1+(a+y)^2]v''(y)=0.
$$
Denote $w(y)=(2a+y)v(y)$. The equation above may be written as
$$
(a+y)w'(y)=[1+(a+y)^2]w''(y).
$$
If $w=$constant then $v(y)=y+2a$.
If $w'(y)\neq0$ then we may solve the ODE obtaining
$$
w(y)=\frac c2(a+y)\sqrt{1+(a+y)^2}+\frac c2\ln(a+y+\sqrt{1+(a+y)^2})+b, \quad b,c\in{\mathbb{R}}.
$$
The case $w=$constant may be included in the previous expression.
Hence
\begin{equation}
\label{sol_v2:eq}
v(y)=\frac b{2a+y}+\frac c{2(2a+y)}\left[(a+y)\sqrt{1+(a+y)^2}+\ln(a+y+\sqrt{1+(a+y)^2})\right].
\end{equation}
We conclude with the following
\begin{theorem}
\label{th2}
Minimal translation surfaces of type $2$ in the Heisenberg group $\mathrm{Nil}_3$ are parameterized by
\begin{equation}
\label{trans2:param}
r(x,y)=(x,0,u(x))*(v(y),y,0)=\Big(x+v(y),y,u(x)+\frac{xy}2\Big),
\end{equation}
where
\begin{itemize}
\item[(i)] either $u(x)=\frac{a}{2(1+a^2)}~x^2+bx+c$ and $v(y)=ay+d$, with $a,b,c,d\in{\mathbb{R}}$,
\item[(ii)] or $u(x)=ax+u_0$ and $v(y)$ is given by \eqref{sol_v2:eq}, with $a, u_0\in{\mathbb{R}}$.
\end{itemize}
\end{theorem}

\begin{remark}\rm
If $a$ and $d$ in (i) vanish (resp. $b$ and $c$ in (ii)), then $M$ is a right cylinder over the curve
$\gamma_1(x)=(x,0,u(x))$. See also \cite{ikos}.
\end{remark}

\subsection{Surfaces of type 5}
Taking the two curves as in the previous case:
$\gamma_1(x)=(x,0,u(x))$ and $\gamma_2(y)=(v(y),y,0)$, the translation surface
$M(\gamma_2,\gamma_1)$ of type 5 is parametrized by
$$
r(x,y)=\gamma_2(y)*\gamma_1(x)=\Big(x+v(y),y,u(x)-\frac{x y}2\Big).
$$
The two parametrizations of $M(\gamma_2,\gamma_1)$ and $M(\gamma_1,\gamma_2)$ look very similar,
yet they are quite different.

Straightforward computations yield the following minimality condition
\begin{equation}
\label{min5:eq}
\begin{array}{ccc}
u''(x)\left[-2y(v(y)+2x)v'(y)+(y^2+4)v'(y)^2+(v(y)+2x)^2+4\right]-4u'(x)^3v''(y)\\[2mm]
+2yu'(x)^2v''(y)-2u'(x)\left[2(v''(y)+x)-yv'(y)+v(y)\right]\\[2mm]
+2\big(yv''(y)+2v'(y)\big)=0.
\end{array}
\end{equation}

Let us make some notations:

\begin{equation}
\label{eq:not_P}
\begin{array}{l}
P_1(x,y)=2x+v(y)-yv'(y),\\[2mm]
P_2(y)=2v'(y),\\[2mm]
P_3(y)=P_2'(y)=2v''(y),\\[2mm]
P_4(y)=2yv''(y).
\end{array}
\end{equation}

The equation\eqref{min5:eq} may be rewritten as

\begin{equation}
\label{min5r:eq}
\left(4+P_1^2+P_2^2\right)u''(x)-2(P_1+P_3)u'(x)+P_4u'(x)^2-2P_3u'(x)^3+(2P_2+P_4)=0.
\end{equation}
Taking the derivative with respect to $x$ we get
\begin{equation}
\label{min5_1:eq}
(4+P_1^2+P_2^2)u'''(x)+2(P_1-P_3)u''(x)-4u'(x)+2P_4u'(x)u''(x)-6P_3u'(x)^2u''(x)=0.
\end{equation}

{\bf (i)} $u''(x)=0$: From \eqref{min5_1:eq} it follows $u'(x)=0$ and hence $u=u_0$, $u_0\in{\mathbb{R}}$.
Combining with \eqref{min5r:eq} we obtain $2P_2+P_4=0$. This ODE has the general solution $v(y)=\frac ay+b$,
$a,b\in{\mathbb{R}}$.

{\bf (ii)} $u''(x)\neq0$: Dividing in \eqref{min5r:eq} by $u''(x)$ and then taking successively the derivatives with respect to
$x$ and $y$ one obtains
$$
2v'(y)v''(y)\frac1{u''(x)}\frac {\rm d}{{\rm d}x}\left(\frac{u'''(x)}{u''(x)}\right)-
yv''(y)\frac{u'''(x)}{u''(x)^2}+\big(v''(y)+yv'''(y)\big)-6v'''(y)u'(x)=0.
$$
Taking one more derivative with respect to $x$ and considering
$$
\begin{array}{l}
\displaystyle
A(x)=\frac 1{u''(x)}\frac{\rm d}{{\rm d}x}\left[\frac1{u''(x)}\frac{\rm d}{{\rm d}x}\left(\frac{u'''(x)}{u''(x)}\right)\right],\\[4mm]
\displaystyle
B(x)=\frac1{u''(x)}\frac{\rm d}{{\rm d}x}\left(\frac{u'''(x)}{u''(x)^2}\right),
\end{array}
$$
we get
$$
2v'(y)v''(y)A(x)-yv''(y)B(x)=6v'''(y).
$$
\begin{enumerate}
\item[(a)] If $v''(y)=0$ then $v(y)=ay+b$, $a,b\in{\mathbb{R}}$. Replacing it in
\eqref{min5r:eq} we obtain the following ODE
$$
[4+4a^2+(2x+b)^2]u''(x)-2(2x+b)u'(x)+4a=0,
$$
with the solution
\begin{equation}
\label{sol_u5}
\begin{array}{rl}
u(x)=&\displaystyle
       -\frac{ax^2+abx+c}{2(1+a^2)}+c_1\left(\frac x2+\frac b4\right)\sqrt{4(1+a^2)+(2x+b)^2}\\[3mm]
     & +c_1(1+a^2)\ln\left(2x+b+\sqrt{4(1+a^2)+(2x+b)^2}\right),\quad c,c_1\in{\mathbb{R}}.
\end{array}
\end{equation}
\item[(b)] If $v''(y)\neq0$, it follows
either $A$ and $B$ are constants,
 or $v'(y)=cy$, $c\neq0$ and $B=2cA$.

\noindent
It is not difficult to show that each of the two situations yields a contradiction.
\end{enumerate}

We conclude with the following
\begin{theorem}
\label{th5}
Minimal translation surfaces of type $5$ in the Heisenberg group $\mathrm{Nil}_3$ are parameterized by
\begin{equation}
\label{trans5:param}
r(x,y)=(v(y),y,0)*(x,0,u(x))=\Big(x+v(y),y,u(x)-\frac{xy}2\Big),
\end{equation}
where
\begin{itemize}
\item[(i)] either $u(x)=u_0$ and $v(y)=\frac ay+b$, with $u_0,a,b\in{\mathbb{R}}$
\item[(ii)] or $v(y)=ay+b$ and $u(x)$ is given by \eqref{sol_u5}.
\end{itemize}
\end{theorem}
\begin{remark}\rm
In the case (i), the curve $\gamma_{1}(x)$ is a geodesic. Thus $M(\gamma_2,\gamma_1)$
is a right translation of $\gamma_2(y)$ by a geodesic $\gamma_{1}(x)$.
If $a$ and $b$ vanish, then the surface $M$ is a left cylinder over the curve $\gamma_1(x)=(x,0,u(x))$.
\end{remark}

\subsection{Surfaces of type 3}
Let $\gamma_1(x)=(0,x,u(x))$ and $\gamma_2(y)=(y,v(y),0)$ be the two curves
defining the translation surface $M(\gamma_1,\gamma_2)$, which is parameterized as
$$
r(x,y)=\gamma_1(y)*\gamma_2(x)=\Big(y,x+v(y),u(x)-\frac{x y}2\Big).
$$
Using the notations given by \eqref{eq:not_T}, the minimality equation may be written as
$$
T_0(y)u''(x)+T_1(y)u'(x)-T_2(y)u'(x)^2+T_3(y)u'(x)^3-T_4(y)=0.
$$
Applying the same technique as in the case of surfaces of type 2, we obtain
\begin{theorem}
\label{th3}
Minimal translation surfaces of type $3$ in the Heisenberg group $\mathrm{Nil}_3$ are parameterized by
\begin{equation}
\label{trans3:param}
r(x,y)=(x,0,u(x))*(v(y),y,0)=\big(x+v(y),y,u(x)+\frac{xy}2\big)
\end{equation}
where
\begin{itemize}
\item[(i)] either $u(x)=\frac{a}{2(1+a^2)}~x^2+bx+c$ and $v(y)=-ay+d$, with $a,b,c,d\in{\mathbb{R}}$,
\item[(ii)] or $u(x)=ax+u_0$ and

\
$v(y)=\displaystyle
\frac b{2a-y}+\frac c{2(2a-y)}\left[(a-y)\sqrt{1+(a-y)^2}+\ln(a-y+\sqrt{1+(a-y)^2})\right],
$
\newline
with $a, b, c, u_0\in{\mathbb{R}}$.
\end{itemize}
\end{theorem}

\subsection{Surfaces of type 6}
Let us consider the two curves as in previous case and define the translation surface
$M(\gamma_2,\gamma_1)=\gamma_2(y)*\gamma_1(x)$.

We can state the following result.
\begin{theorem}
\label{th6}
Minimal translation surfaces of type $6$ in the Heisenberg group $\mathrm{Nil}_3$ are parameterized by
\begin{equation}
\label{trans6:param}
r(x,y)=(y,v(y),0)*(0,x,u(x))=\Big(y,x+v(y),u(x)+\frac{xy}2\Big),
\end{equation}
where
\begin{itemize}
\item[(i)] either $u(x)=u_0$ and $v(y)=\frac ay+b$, with $u_0,a,b\in{\mathbb{R}}$
\item[(ii)] or $v(y)=ay+b$ and

$\begin{array}{rl}
u(x)=& \displaystyle
\frac{ax^2+abx+c}{2(1+a^2)}+c_1\left(\frac x2+\frac b4\right)\sqrt{4(1+a^2)+(2x+b)^2}\\[2mm]
& +c_1(1+a^2)\ln\left(2x+b+\sqrt{4(1+a^2)+(2x+b)^2}\right),
\end{array}
$
\newline
with $a, b,c,c_1\in{\mathbb{R}}$.
\end{itemize}
\end{theorem}

\begin{remark}\rm
We may construct minimal translation surfaces which are cylinders also for types 3 and 6.
\end{remark}

\begin{remark}\rm
Minimal translation surfaces of types 1-4 are also flat.
\end{remark}

\begin{remark}\rm Some missing cases:

In general, when we consider a curve in a certain plane of coordinates, parameterized by $s$
(not necessary the arc-length), namely $s\longmapsto \big(\alpha(s),\beta(s)\big)$, with
$\alpha'(s)^2+\beta'(s)\neq0$, we may write the curve also in the explicit form, as
$(t,u(t))$ and this can be done if $\alpha'(s)\neq0$. So, we have also to consider curves given by
$(c,t)$, corresponding to $\alpha'(s)=0$ and $\beta'(s)\neq0$.

Having this in mind, we describe the missing cases in our classifications:


\begin{table}[htb]
\centering
\begin{tabular}{ccccc}
  \toprule
  type & parametrization && type & parametrization \\\midrule[\heavyrulewidth]

  1\ & $\big(c,y,x+v(y)+\frac{c y}2\big)$ && 4\ & $\big(c,y,x+v(y)-\frac{c y}2\big)$ \\
  \midrule[0.1mm]
  1\ & $\big(x,c,u(x)+y+\frac{c x}2\big)$ && 4\ & $\big(x,c,u(x)+y-\frac{c x}2\big)$ \\
  \midrule[0.1mm]

  2$^*$ & $\big(c+v(y),y,x+\frac{c y}2\big)$ && 5$^*$ & $\big(c+v(y),y,x-\frac{c y}2\big)$ \\
  \midrule[0.1mm]
  2 & $\big(x+y,c,u(x)+\frac{c x}2\big)$ && 5 & $\big(x+y,c,u(x)-\frac{c x}2\big)$ \\
  \midrule[0.1mm]

  3$^*$ & $\big(y,c+v(y),x-\frac{c y}2\big)$ && 6$^*$ & $\big(y,c+v(y),x+\frac{c y}2\big)$ \\
  \midrule[0.1mm]
  3 & $\big(c,x+y,u(x)-\frac{c x}2\big)$ && 6 & $\big(c,x+y,u(x)+\frac{c x}2\big)$ \\
  \bottomrule
\end{tabular}
\end{table}

Except the cases marked by $*$ symbol, the others are all minimal translation surfaces.
The stared cases are minimal if and only if $v$ is an affine function.
\end{remark}

{\bf Acknowledgement.}
The first author was partially supported by Kakenhi 21546067.
The second author was partially supported by MEC-FEDER grant no. MTM2007-61775 and
   Junta de Andaluc\'{\i}a grant no. P06-FQM-01642.
The last author was supported by Fulbright Grant n. 498/2010 as a Fulbright Senior Researcher at the Michigan State University, USA.
He wishes to thank Prof. B.Y. Chen for kind hospitality and encouragement during his stay at MSU.



\begin{thebibliography}{99}

\bibitem{ar}
U.~Abresch, H.~Rosenberg,
\emph{Generalized Hopf differentials},
Mat. Contemp. \textbf{28} (2005), 1--28.

\bibitem{bt}
D.~A.~Berdinski{\u\i},
I.~A.~Ta{\u\i}manov,
\emph{Surfaces in three-dimensional
Lie groups}, Siberian Math. J.
\textbf{46} (2005) 6, 1005--1019.

\bibitem{d} B.~Daniel,
\emph{The Gauss map of minimal surfaces in the Heisenberg group},
Int Math Res Notices (2011) 3, 674--695.


\bibitem{dhm} B.~Daniel, L.~Hauswirth, P.~Mira,
\emph{Constant mean curvature surfaces in homogeneous $3$-manifolds},
Lectures Notes of the $4^{\mathrm{th}}$ KIAS Workshop on Diff. Geom.
{\sl Constant mean curvature surfaces in homogeneous manifolds},
Seoul, 2009.

\bibitem{dgv}
F.~Dillen, W.~Goemans, I.~Van de Woestyne,
\emph{Translation surfaces of Weingarten type in 3-space},
{\sl Proc. Conf. RIGA 2008},
Bull. Transilvania Univ. Brasov, \textbf{15} (2008) 50,
1--14.

\bibitem{f} C.~B.~Figueroa,
\emph{On the Gauss map of a minimal surface in the Heisenberg group},
Mat. Contemp. \textbf{33} (2007),
139--156.

\bibitem{fm} I.~Fern{\'a}ndez, P.~Mira,
\emph{Constant mean curvature surfaces in
$3$-dimensional Thurston geometries},
to appear in Proceedings on the ICM 2010,
Hyderabad, arXiv:1004.4752v1 [math.DG] (2010).

\bibitem{gv}
W.~Goemans, I.~Van de Woestyne,
\emph{Translation surfaces with vanishing second Gaussian curvature
in Euclidean and Minkowski $3$-space},
Proceedings of the conference
Pure and Applied Differential Geometry,
PADGE 2007, Eds. F. Dillen, I. Van de Woestyne, 123--131.

\bibitem{i2}
J.~Inoguchi,
\emph{Flat translation invariant surfaces
in the $3$-dimensional Heisenberg group},
J. Geom. \textbf{82} (2005), 83--90.

\bibitem{i3} J.~Inoguchi,
 \emph{Minimal surfaces in
the $3$-dimensional Heisenberg group},
Diff. Geom. Dyn. Syst.
\textbf{10} (2008), 163--169.

\bibitem{ikos}
J.~Inoguchi, T.~Kumamoto,
N.~Ohsugi, Y.~Suyama,
\emph{Differential geometry of curves
and surfaces in $3$-dimensional homogeneous spaces}, II,
Fukuoka Univ. Sci. Reports \textbf{30} (2000) 1, 17--47.

\bibitem{im}
J.~Inoguchi, M.~I.~Munteanu,
\emph{Minimal translation surfaces in the hyperbolic $3$-space},
in preparation.

\bibitem{lm}
R.~L{\'o}pez, M.~I.~Munteanu,
\emph{Minimal translation surfaces in $\mathrm{Sol}_3$},
to appear in J. Math. Soc. Japan.

\bibitem{mn}
M.~I.~Munteanu, A.~I.~Nistor,
\emph{On the geometry of the second fundamental form of translation
surfaces in $\mathbb{E}^3$}, Houston J. Math. \textbf{37} (2011) 4, xxx - xxx.

\bibitem{sc} H.~F.~Scherk,
\emph{Bemerkungen \"{u}ber die kleinste Fl\"{a}che innerhalb gegebener Grenzen},
J. Reine Angew. Math. \textbf{13} (1835), 185--208.




\bibitem{th}  W.~Thurston,
\emph{Three-dimensional geometry and topology}, Princeton Math. Ser. 35,
Princeton Univ. Press, Princeton, NJ, 1997.




\end{thebibliography}
\end{document}